%

\magnification=1200
\input amstex
\documentstyle{amsppt}


\topmatter
\title
An Induction Principle and\\
Pigeonhole Principles for K-finite Sets 
\endtitle
\rightheadtext{Induction and Pigeonhole Principles}
\author
Andreas Blass
\endauthor
\address
Mathematics Dept., University of Michigan, Ann Arbor, 
MI 48109, U.S.A.
\endaddress
\email
ablass\@umich.edu
\endemail
\thanks Partially supported by NSF grant DMS-9204276.
\endthanks
\subjclass
03F55
\endsubjclass
\abstract
We establish a course-of-values induction principle for K-finite sets
in intuitionistic type theory.  Using this principle, we prove a
pigeonhole principle conjectured by B\'enabou and Loiseau.  We also
comment on some variants of this pigeonhole principle.
\endabstract
\endtopmatter
\document

\head
1. Introduction
\endhead

The pigeonhole principle says that a finite set cannot be mapped
one-to-one into a proper subset.  There is a dual principle saying
that a finite set cannot be mapped onto a proper superset.  We
consider these principles in the context of constructive logic.  The
motivation for these considerations came from a weak version of the
dual pigeonhole principle proved constructively by B\'enabou and
Loiseau, who noted that their argument does not establish a natural
stronger version of the principle.

Throughout this paper we work in an intuitionistic type theory of the
sort that arises as the internal logic of an elementary topos \cite{3,
5, 6, 7, 8, 9, 11}.  Although the questions we consider originated in
the course of topos-theoretic work of B\'enabou and Loiseau \cite4,
many of our results are theorems of intuitionistic type theory (in
fact of intuitionistic third-order logic) and involve no reference to
topoi.

Of the several concepts of finiteness that are equivalent in classical
logic but not in intuitionistic logic, we shall use the one commonly
called K-finiteness or Kuratowski-finiteness \cite{1, 9, 10}.  The
definition and some comments on it are given in Section~2.
Henceforth, we omit the prefix K and refer simply to finiteness.

B\'enabou and Loiseau showed \cite{4, Prop.~5.3} that, if $X$ is
finite and inhabited, then no function $f:X\to X\times2$ can be
surjective.  In other words, for every such $f$ it is not the case
that every element of $X\times2$ is in its range.  They pointed out
that this version of the dual pigeonhole principle is weaker (in
constructive logic) than the statement that for every such $f$ there
is an element of $X\times2$ not in its range.  They remarked that the
latter, stronger statement ``seems to be true, but we do not have a
general proof of it.''  One purpose of the present paper is to give a
general proof of it and in fact of the stronger statement obtained by
replacing $X\times2$ with $X+1$.

The proof uses an induction principle whereby, when one proves a
property for an arbitrary finite set $X$, one can assume the property
for all complemented, proper subsets of $X$.  This principle seems to
be of interest independently of the application that motivated it.

After explaining our terminology and presenting some preliminary facts
in Section~2, we devote Section~3 to proving the induction principle.
The application to the stronger version of the dual pigeonhole
principle conjectured by B\'enabou and Loiseau is in Section~4.  In
Section~5, we consider the (undualized) pigeonhole principle, showing
that a weak version is intuitionistically provable but a strong
version is not.  The final Section~6 is about some variants of the
dual pigeonhole principle.

\head
2. Preliminaries
\endhead

The set theory and logic used in this paper are an intuitionistic type
theory of the sort described in \cite{3, 5, 6, 7, 8, 9, 11}.  These
references differ in some details, but the differences will not matter
in our work.  We shall work with elements, subsets, and families of
subsets of some fixed but arbitrary type $U$, as well as (partial)
functions from $U$ to $U$.

A set is called {\sl inhabited\/} if there exists an element in it.
This is stronger (in intuitionistic logic) than not being empty.  A
subset $A$ of a set $B$ is a {\sl proper\/} subset if $B-A$ is
inhabited, i.e., if $B$ has an element that is not in $A$.  We say
that $A$ is {\sl complemented\/} in its superset $B$ if
$B=A\cup(B-A)$, i.e., if every element of $B$ is either in $A$ or not
in $A$.  

A function $f$ is {\sl one-to-one\/} if $f(x)=f(y)$ implies $x=y$.
This definition is intuitionistically stronger and more natural than
the classically equivalent but negation-filled definition that $x\neq
y$ implies $f(x)\neq f(y)$.

We call a set $A$ {\sl finite\/} if it belongs to every family $\Cal
X$ that contains the empty set $\emptyset$ and is closed under
adjoining single elements in the sense that if $Z\in\Cal X$ and $a\in
A$ then $Z\cup\{a\}\in\Cal X$.  

The definition trivially implies that $\emptyset$ is finite and that,
if $A$ is finite, then so is $A\cup\{p\}$ for every $p$.

This definition also immediately implies an induction principle.  To
prove that all finite sets have some property, it suffices to prove
that the empty set has the property (induction basis) and that,
whenever $Z$ has the property and $a$ is any element (in $U$), then
$Z\cup\{a\}$ also has the property (induction step).  We shall refer
to this sort of induction as {\sl ordinary\/} induction on finite
sets, to distinguish it from the new induction principle to be
established in Section~3.

The definition of finiteness is (intuitionistically) equivalent (cf.
\cite{4, Lemma~5.2}) to the definition of K-finiteness given in
\cite9, namely that $A$ belongs to every family $\Cal Y$ that contains
$\emptyset$ and $\{a\}$ for all $a\in A$ and is closed under binary
union.  To see the equivalence, note first that any $\Cal Y$ as in the
second definition is also an $\Cal X$ as in the first definition, so
all finite sets are K-finite.  For the other direction, one shows that
the family of finite sets is closed under binary union, so the finite
subsets of $A$ form a $\Cal Y$ as in the second definition.  To prove
that, if $x$ and $y$ are finite then so is $x\cup y$, one proceeds by
ordinary induction on $x$; both the basis and the induction step are
trivial.

An equally trivial induction establishes that every finite set is
either empty or inhabited.  

A complemented subset $A$ of a finite set $B$ is finite.  To see this,
proceed by induction on $B$, the basis ($B=\emptyset$) being trivial.
So suppose the result is true for $B$ and that $A$ is a complemented
subset of $B\cup\{p\}$.  Then $A\cap B$ is a complemented subset of
$B$, so by induction hypothesis it is finite.  If $p\notin A$, then
$A=A\cap B$ and thus $A$ is finite.  If $p\in A$ then $A=(A\cap
B)\cup\{p\}$ and again $A$ is finite.  Since $A$ is complemented in
$B\cup\{p\}$, the cases considered in the preceding two sentences
exhaust the possibilities, so the proof is complete.

It follows from the preceding two paragraphs that, if $A$ is a
complemented subset of $B$ then either it is a proper subset or it
equals $B$.  Indeed, $B-A$ is also complemented, hence finite, and
hence either inhabited or empty.  If $B-A$ is inhabited, then $A$ is a
proper subset of $B$.  If $B-A$ is empty then, as $A$ is complemented
in $B$, we have $B=A\cup(B-A)=A$. 

We emphasize that the finite sets we work with need not have a
decidable equality relation.  That is, we do not assume that $x=y$ or
$x\neq y$.  In fact, by \cite1, such an assumption would allow us to
work in a sub-universe (sub-topos) in which classical logic holds and
would thus remove the whole point of working in intuitionistic logic.

\head
3. An Induction Principle
\endhead

This section is devoted to establishing an induction principle,
different from the one given by the definition of finiteness, for
proving properties of finite sets.  

\proclaim{Theorem 1}
Let $\Cal X$ be a family of finite sets such that 
\roster
\item for all finite $A$, if all complemented proper subsets of $A$ 
are in $\Cal X$, then also $A\in\Cal X$.
\endroster
Then $\Cal X$ contains all finite sets.
\endproclaim

The theorem says that, in order to prove a statement for all finite
sets, it suffices to prove it for an arbitrary finite set $A$ assuming
that it holds for all complemented proper subsets of $A$.  It is
related to ordinary induction for finite sets much as course-of-values
induction is related to ordinary induction for natural numbers.  

\demo{Proof}
We show, by ordinary induction on finite sets $B$ (in the sense
explained in Section~2) that 
$$
\forall\Cal X\,[\text{\therosteritem1}\implies B\in\Cal X].\tag2
$$

The basis is easy, for if $B=\emptyset$ then $B$ has no proper subset,
so \therosteritem1 applied with $A=B$ immediately gives $B\in\Cal X$.

For the induction step, we assume \thetag2 for a particular finite
$B$; we wish to prove \thetag2 for $B\cup\{p\}$.  So fix an $\Cal X$
satisfying \therosteritem1; we must prove $B\cup\{p\}\in\Cal X$.
Define 
$$
\Cal Y=\{A\mid A\in\Cal X \text{ and } A\cup\{p\}\in\Cal X\}.
$$

We claim that \therosteritem1 holds with $\Cal Y$ in place of $\Cal
X$.  Assuming the claim for a moment, we can apply the induction
hypothesis \thetag2 for $B$ with $\Cal X$ instantiated as $\Cal Y$.
So we get $B\in\Cal Y$, from which the desired $B\cup\{p\}\in\Cal X$
immediately follows by definition of $\Cal Y$.

So all that remains is to prove the claim that \therosteritem1 holds
with $\Cal Y$ in place of $\Cal X$.  So let $A$ be finite and assume
that all its complemented proper subsets are in $\Cal Y$.  In
particular, all its complemented proper subsets are in $\Cal X$ and so
$A\in\Cal X$ since \therosteritem1 holds for $\Cal X$.  It remains to
prove that $A\cup\{p\}\in\Cal X$, and we shall do this by applying the
assumption \therosteritem1 for $\Cal X$.

So let $C$ be any complemented proper subset of $A\cup\{p\}$; we must
show $C\in\Cal X$.  As $C$ is complemented, we have $p\in C$ or
$p\notin C$.  Also, as $C\cap A$ is a complemented subset of $A$, it
is either equal to $A$ or a proper subset of $A$.  We consider the
various cases.

If $C\cap A$ is a proper subset of $A$, we use the assumption that
$\Cal Y$ contains all the complemented proper subsets of $A$ to
conclude that $C\cap A\in\Cal Y$.  Then $C$, being equal to ($C\cap
A)\cup\{p\}$ or to $C\cap A$ (according to whether $p\in C$), is in
$\Cal X$ by definition of $\Cal Y$.

If $C\cap A=A$ and $p\notin C$ then $C=A$, and we already saw that
$A\in\Cal X$.

The remaining case, $C\cap A=A$ and $p\in C$, is impossible as $C$ is a
proper subset of $A\cup\{p\}$.

Thus, we have $C\in\Cal X$ in all cases, which completes the proof.
\qed\enddemo

Although Theorem~1 suffices for the proofs in the following sections,
it seems natural to ask whether it could be strengthened by replacing
``complemented'' with ``finite'' in \therosteritem1.  It is not
difficult to prove this strengthened induction principle if the axiom
of infinity is available, that is, if the type $\Bbb N$ of natural
numbers is available in the intuitionistic type theory.  The proof
begins by showing, by ordinary induction on finite sets $A$, that
there is a natural number $n$ such that $\bar n=\{0,1,\dots,n-1\}$
has no one-to-one map into $A$.  Then one shows by induction on $n$
that any set $A$ admitting no one-to-one map from $\bar n$ must be in
every class $\Cal X$ that satisfies the weakened version of
\therosteritem1.  

It is not clear to me whether one can obtain the same result without
an axiom of infinity, but it seems that any proof would have to be
substantially different from the one just given.  To see this,
consider the statement ``If $A$ is finite then there is a finite $B$
such that the equality relation on $B$ is decidable and $B$ has no
one-to-one map into $A$.''  This statement is a reformulation, in the
absence of $\Bbb N$, of the result of the first half of the proof
given above.  (See \cite1 for the connection between the sets $\bar n$
and finite sets $B$ with decidable equality.)  But this statement is
not provable in intuitionistic type theory without the axiom of
infinity.  More precisely, there exist a topos $\Cal E$ and a finite
object $A$ in it such that, for any object $C$ of $\Cal E$, the
statement ``$C$ has a finite subset with decidable equality admitting
no one-to-one map into $A$'' fails to be internally valid.  

To construct such a topos, let $P$ be a three-element partially
ordered set with a top element 1 and two incomparable elements $a$ and
$b$ below it.  Let $\Cal E'$ be the topos of presheaves on $P$ in some
non-standard model of set theory, and let $\Cal E$ be the subtopos
consisting of those presheaves whose values at 1 and $a$ are finite in
the sense of that non-standard model and whose values at $b$ are
really finite.  Fix a set $S$ that is finite in the sense of the
non-standard model but is not really finite, and let $A$ be the
presheaf whose values at 1 and $a$ are $S$ with the identity as
transition map between them and whose value at $b$ is a singleton.
This $A$ is finite in $\Cal E$.  If $C$ is any other object of $E$,
then one can calculate, using Kripke-Joyal semantics, that a finite
subset of $C$ with decidable equality at 1 would have all three of its
components really finite (the $b$ component by definition of $\Cal E$,
then the 1 component because decidability makes the transition maps
one-to-one, and then the $a$ component because finiteness makes the
transition maps surjective).  So, after restriction to $a$, it could
be mapped one-to-one into $A$ since $A(a)$ is really infinite.

\head
4. The Strong Dual Pigeonhole Principle
\endhead

The purpose of this section is to establish, in intuitionistic type
theory, the stronger version of the dual pigeonhole principle
conjectured by B\'enabou and Loiseau, namely that if $X$ is finite and
inhabited and $f:X\to X\times2$ then there is an element of $X\times2$
that is not in the range of $f$.  In fact, our proof gives a stronger
statement with $X+1$ instead of $X\times2$.  (To see that the $X+1$
result is indeed stronger than the $X\times2$ result, it suffices to
observe that, since $X$ is inhabited, there is a surjection
$X\times2\to X+1$ sending one copy of $X$ in $X\times2$ onto $X$ and
the other copy onto 1.)

\proclaim{Theorem 2}
If $A$ is a finite, complemented, proper subset of $B$ and if $f:A\to
B$, then $B-\text{Range}(f)$ is inhabited.
\endproclaim

\demo{Proof}
Let $\Cal X$ be the family of those finite sets $A$ such that, for
every $B$ in which $A$ is a complemented, proper subset and for every
$f:A\to B$, there is an element of $B$ not in the range of $f$.  We
prove that $\Cal X$ contains all finite sets by applying Theorem~1.
So it suffices to prove $A\in\Cal X$ under the assumptions that $A$ is
finite and that every complemented, proper subset of $A$ is in $\Cal
X$.

To do this, suppose $f:A\to B$ where $A$ is a complemented, proper
subset of $B$.  Since $A$ is complemented in $B$, $f^{-1}(A)$ is
complemented in $f^{-1}(B)=A$.  As was pointed out in Section~2, it
follows that $f^{-1}(A)$ either equals $A$ or is a proper subset of
$A$.  If $f^{-1}(A)=A$, then $\text{Range}(f)\subseteq A$, so the
complement in $B$ of this range includes $B-A$, which is inhabited
because $A$ is a proper subset of $B$.  So the desired conclusion
holds in this case.

There remains the case that $f^{-1}(A)$ is a complemented, proper
subset of $A$ and is therefore in $\Cal X$.  Apply the definition of
$\Cal X$ with $f^{-1}(A)$, $A$ and $f\restriction f^{-1}(A)$ in the
roles of $A$, $B$, and $f$.  It shows that
$A-\text{Range}(f\restriction f^{-1}(A))$ is inhabited.  But this set
equals $A-\text{Range}(f)\subseteq B-\text{Range}(f)$.
\qed\enddemo

\head
5. The Undualized Pigeonhole Principle
\endhead

In this section we consider the principle that a finite set cannot be
mapped one-to-one into a proper subset.  More precisely, we consider
two intuitionistically inequivalent versions of this principle.  The
weaker version says that if $X$ is finite then a map $f:X+1\to X$
cannot be one-to-one.  The stronger version says that if $X$ is
finite and $f:X+1\to X$ then there exist $x,y\in X+1$ such that $x\neq
y$ but $f(x)=f(y)$.  

Of course in classical logic these are equivalent and easy to prove.
We shall show that the weaker version is intuitionistically provable
but the stronger is not.  In fact, the stronger version implies the
law of the excluded middle.

\proclaim{Theorem 3}
If $X$ is finite then there is no one-to-one function from $X+1$ into
$X$. 
\endproclaim

\demo{Proof}
The statements ``$X$ is finite,'' ``$f:X+1\to X$,'' and ``$f$ is
one-to-one'' are all preserved by inverse images of geometric
morphisms of topoi (see \cite{9} especially Cor.~9.17).  So if their
conjunction were intuitionistically consistent and therefore had a
non-zero truth value in some elementary topos, then, by Barr's theorem
\cite{2, 9}, it would have a non-zero truth value in some Boolean topos.
That is absurd, since the pigeonhole principle is provable in
classical type theory and therefore valid in every Boolean topos.  So
the conjunction of the three statements is intuitionistically
inconsistent.
\qed\enddemo

We remark that a similar proof can be given for the weak form of the
dual pigeonhole principle.  The statements ``$X$
is finite,'' ``$f:X\to X+1$,'' and ``$f$ is surjective'' are preserved
by inverse images of geometric morphisms, so Barr's theorem allows us
to conclude their intuitionistic inconsistency from their classical
inconsistency.

In fact, we can do a bit better and replace ``$X$ is finite'' by the
intuitionistically weaker ``$X$ is a subset of a finite set'' and
still conclude that there is no surjection $X\to X+1$ and no
one-to-one map $X+1\to X$.  This is because the weaker hypothesis
suffices for the classical proof and its internal validity is
preserved by geometric inverse images.

In contrast to the situation with the dual pigeonhole principle, where
the stronger form turned out to be provable (Theorem~2), the
undualized pigeonhole principle cannot be similarly strengthened
without going to classical logic.

\proclaim{Theorem 4}
Assume that, for all finite $X$ and all $f:X+1\to X$, there exist $x$
and $y$ in $X+1$ with $f(x)=f(y)$ but $x\neq y$.  Then the law of the
excluded middle holds.
\endproclaim

\demo{Proof}
Let $u$ be an arbitrary truth value, and let $X$ be a set whose
elements are exactly $a$ and $b$ where $a=b$ if and only if $u$.  Such
a set $X$ can be obtained as the quotient of $1+1$ by an equivalence
relation containing all pairs if $u$ and also containing the diagonal
pairs; the equivalence classes of the two distinct elements of $1+1$
serve as $a$ and $b$.  Clearly, $X$ is finite.

Writing $c$ for the unique element of 1, we define a map $f$ from
$X+1=\{a,b,c\}$ to $X=\{a,b\}$ by sending $a$ to itself, $b$ to $a$,
and $c$ to $b$.  This is well-defined even though $a$ might equal $b$,
since they are sent to the same element $a$.  

By assumption, there are $x,y\in\{a,b,c\}$ with $f(x)=f(y)$ but $x\neq
y$.  We have $x=a$ or $x=b$ or $x=c$ and similarly for $y$, so we can
consider the nine resulting (exhaustive though not necesarily
exclusive) cases.  Three ``diagonal'' cases have $x=y$ contrary to the
choice of $x$ and $y$.  Two other cases have $x=c$ while $y=a$ or
$y=b$; in these cases $f(x)=f(y)$ means that $b=a$ and therefore $u$
holds.  The two similar cases with $y=c$ also give that $u$ holds.
There remain two cases, one with $x=a$ and $y=b$ and the symmetric one
with $x=b$ and $y=a$.  In either of these two cases, $x\neq y$ means
that $a\neq b$ and therefore not $u$.  Thus, in all cases, we have $u$
or not $u$.  As $u$ was an arbitrary truth value, the proof is
complete.
\qed\enddemo

\head
6. Variants of the Dual Pigeonhole Principle
\endhead

This section is devoted to refuting two possible strengthenings of
Theorem~2.  The first is to weaken the hypothesis from ``finite'' to
``subset of a finite set,'' as we did with the weak pigeonhole
principles in the remarks following Theorem~3.  For the strong dual
pigeonhole principle, this further strengthening is not only
unprovable intuitionistically but equivalent to classical logic.

\proclaim{Theorem 5}
Assume that, whenever $A$ is a subset of a finite set and $f:A\to A+1$
then there is an element of $A+1$ not in the range of $f$.  Then the
law of the excluded middle holds.
\endproclaim

\demo{Proof}
Let $u$ be any truth value, and let $U$ be the corresponding subobject
of 1 (inhabited if and only if $u$).  Since 1 is finite, we can apply
the hypothesis of the theorem with $A=U$.  Let $f:U\to U+1$ be the
inclusion of $U$ in the {\sl second\/} summand 1 of $U+1$.  By
hypothesis, $U+1$ has an element $x\notin\text{Range}(f)$.  If $x$ is
in the first summand $U$, then (as $U$ is thereby inhabited) $u$
holds.  If $x$ is in the second summand 1, then, by definition of $f$,
it would belong to the range of $f$ with truth value $u$.  Since it
does not belong to this range, we conclude that not $u$.  As $x$ must
be in one of the two summands, we have proved that $u$ or not $u$.
\qed\enddemo

Finally, we consider an attempt to extend Theorem~2 from the internal
logic of topoi (intuitionistic type theory) to the external logic.
Specifically, if $A$ is a finite object in a topos $\Cal E$ and if
$f:A\to A+1$ is a morphism in $\Cal E$, must $A+1$ have a global
section disjoint from the image of $f$?  The answer is negative.  For
a counterexample, consider a topos with an inhabited finite object $A$
having no global section, and let $f:A\to A+1$ map $A$ onto the
summand 1.  A global section of $A+1$ disjoint from the range of $f$
would be a global section of $A$ and thus does not exist.  For a
simple example of a topos containing such an $A$, use the topos of
presheaves on a four-element poset $\{p,q,r,s\}$ where $p$ is
incomparable with $q$, $r$ is incomparable with $s$, and both $p$ and
$q$ are below both $r$ and $s$.  Let $A$ be the presheaf whose value
at each point is $1+1$ and whose transition maps are the identity map
of $1+1$ except for one, say from $r$ to $p$, that interchanges the
elements of $1+1$.  Then $A$ clearly has no global section, but it is
inhabited and finite and in fact internally isomorphic to $1+1$.

\enddocument